\documentclass[a4paper,10pt]{amsart}
\usepackage[utf8]{inputenc}
\usepackage{graphicx}
\usepackage[american]{babel}
\usepackage{wasysym}

\usepackage{geometry}
\usepackage{amsmath}
\usepackage{amssymb}
\usepackage{latexsym}

\usepackage{amsthm}

\usepackage{indentfirst}
\usepackage{color}

\usepackage{ifpdf}

\usepackage{framed}
\usepackage{bookmark}
\usepackage{graphicx}
\usepackage{amsfonts}
\usepackage{graphicx}
\usepackage{amsfonts}
\usepackage{bbm}
\usepackage{booktabs}
\usepackage{times}
\usepackage{float}
\usepackage[all]{xy}
\usepackage{epstopdf}
\usepackage{mathrsfs}
\usepackage{prettyref}
\theoremstyle{plain}

\theoremstyle{definition}\newtheorem{definition}{Definition}[section]

\newrefformat{lemma}{Lemma \ref{#1}}
\newrefformat{thm}{Theorem \ref{#1}}
\newrefformat{prop}{Proposition \ref{#1}}
\newrefformat{rem}{Remark \ref{#1}}
\newrefformat{exercise}{Exercise \ref{#1}}

\ifpdf
\DeclareGraphicsExtensions{.pdf, .jpg, .tif, .mps}
\else
\DeclareGraphicsExtensions{.mps, .jpg, .eps}
\fi

\title[Non-existence of an Algorithm for Determining the Triviality of $\pi_2$]{Non-existence of an Algorithm for Determining the Triviality of the Second Homotopy Groups of Arbitrary Compact Topological Spaces}
\author{Lefit Yuxiang Hao\textsuperscript{\rm 1}}\address{\textsuperscript{\rm 1}College of Mathematics Science, Tongji University, Shanghai 200092, China}\email{haoyuxiang@tongji.edu.cn}

\author{Zijie Kang\textsuperscript{\rm 2}}\address{\textsuperscript{\rm 2}School of Mathematics and Statistics,Wuhan University of Technology,Wuhan 430070,China}\email{327334@whut.edu.cn}

\author{Hongjie Liu\textsuperscript{\rm 3}}\address{\textsuperscript{\rm 3}Faculty of Science and Technology, BNU-HKBU United International College, Zhuhai 519000, China}\email{r130033019@mail.uic.edu.cn}

\author{Pengcheng Ma\textsuperscript{\rm 4}}\address{\textsuperscript{\rm 4}Shenyuan Honors College, Beihang University, Beijing 100191, China}\email{23230504mpc@buaa.edu.cn}

\author{Mufeng Zhou\textsuperscript{\rm 5}}\address{\textsuperscript{\rm 5}Department of Mathematical science, the university of Nottingham Ningbo China, Ningbo 315199, China}\email{smymz3@nottingham.edu.cn}
\begin{document}
	\begin{abstract}
		In this paper, we demonstrate the non-existence of a computational algorithm capable of determining whether the second homotopy group of any compact constructive topological space is trivial. This finding shows the inherent limitations of computer methods in resolving certain topological problems.
	\end{abstract}
	\maketitle
	
	\section{Introduction}
	The study of homotopy groups is fundamental in algebraic topology, providing deep insights into the structure and classification of topological spaces. In particular, the second homotopy group $\pi_2(X)$ of a topological space captures essential information about the 2-dimensional spheres in this space. While determining the triviality of $\pi_2(X)$ is a significant problem, it remains an elusive challenge, especially for constructive topological spaces. For such a problem, constructive mathematics seems ideal for providing methods of exploring related issues.

	The field of constructive mathematics was developed by two different schools. The Russian school, led by Markov \cite{markov1, markov2} and Shanin \cite{shanin}, and detailed by Kushner \cite{kushner}, incorporates Markov’s principle. This principle states that if assuming a decidable subset of natural numbers is empty leads to a contradiction, one can then identify an element of this set, which is a broader approach compared to Bishop's constructivism. On the other hand, E. Bishop \cite{bishop} and his followers developed a related but distinct approach to constructive mathematics.
	
	\section{Definitions}

	\begin{definition}[Word]
		Given a set $\Sigma$ as the alphabet, a \textit{word} is formed by concatenating a finite number of elements. For example, “0010110” is a word of alphabet $\{0, 1\}$.
	\end{definition}
\begin{definition}[Decidable]
	Let $X$ be a set of some words of alphabet $\Sigma$. $X$ is called decidable if there exists an algorithm that terminates on every word $w$ of $\Sigma$ and gives answer “YES” exactly if $w \in X$ and it gives answer “NO” otherwise.
	
	So a set $X$ of natural numbers is called decidable if there exists an algorithm that determines whether an arbitrarily given natural number $n$ belongs to the set $X$. It is proved that there always exist undecidable sets \cite{shen}.
	
\end{definition}
\begin{definition}[$\varepsilon$-net]
	For a compact metric space $(X, d)$ and a given real number $\varepsilon > 0$, an $\varepsilon$-net is a discrete subset $N \subseteq X$ such that for every point $x \in X$, there exists a point $y \in N$ with $d(x, y) < \varepsilon$. The existence of finite $\varepsilon$-net for all $\varepsilon$ is given by the compactness of the space. The existence of an algorithmically given finite $\varepsilon$-net for any rational $\varepsilon$ can be used as a definition of a constructive compact topological space.
\end{definition}
\begin{definition}[Group]
	A \textit{group} is a nonempty set $S$ with a closed binary operation $\cdot$ such that:
	\begin{itemize}
		\item $(a \cdot b) \cdot c = a \cdot (b \cdot c) \quad \forall a, b, c \in S$
		\item $\exists e \in S$ such that $\forall a \in S, a \cdot e = e \cdot a = a$
		\item $\forall a \in S, \exists b \in S$ such that $a \cdot b = b \cdot a = e$
	\end{itemize}
	A group $S$ is trivial if $S$ consists only of the identity element, i.e., $S = \{e\}$.
\end{definition}
\begin{definition}[$n$-cube]
	Let $I^n (n \geq 1)$ denote the unit $n$-cube $I \times \cdots \times I$:
	\[
	I^n = \{(s_1, \ldots, s_n) \mid 0 \leq s_i \leq 1, 1 \leq i \leq n\}
	\]
	The boundary $\partial I^n$ is defined as the geometrical boundary of $I^n$:
	\[
	\partial I^n = \{(s_1, \ldots, s_n) \in I^n \mid \exists i, s_i = 0 \text{ or } 1\}
	\]
\end{definition}
\begin{definition}[$n$-sphere]
	Let $\alpha : I^n \to X$ be a continuous map. $\alpha$ is called $n$-sphere at $x_0$ if $\alpha$ maps the boundary $\partial I^n$ to a point $x_0 \in X$.
\end{definition}
	\begin{definition}[Homotopy]
		Let $X$ be a topological space and $\alpha, \beta : I^n \to X$ be $n$-spheres at $x_0 \in X$. The map $\alpha$ is homotopic with $\beta$, denoted by $\alpha \simeq \beta$, if there exists a continuous map $F : I^n \times I \to X$ such that
	$$\begin{aligned}
			F(s_1, \ldots, s_n, 0) &= \alpha(s_1, \ldots, s_n) \\
			F(s_1, \ldots, s_n, 1) &= \beta(s_1, \ldots, s_n) \\
			F(s_1, \ldots, s_n, t) &= x_0 \quad \text{for } (s_1, \ldots, s_n) \in \partial I^n, t \in I
		\end{aligned}$$
		$F$ is called a homotopy between $\alpha$ and $\beta$. It is easy to show that $\alpha \simeq \beta$ is an equivalence relation. The equivalence class to which $\alpha$ belongs is called the homotopy class of $\alpha$ and is denoted by $[\alpha]$.
	\end{definition}
\begin{definition}[Homotopy Group]
		Let $X$ be a topological space. The set of homotopy classes of $n$-spheres at $x_0 \in X$ is denoted by $\pi_n(X, x_0)$ and called the $n$th homotopy group at $x_0$. The product of homotopy classes $[\alpha]$ and $[\beta]$ is defined by $\alpha \ast \beta = [\alpha \cdot \beta]$. The product of two $n$-spheres $f$ and $g$ is defined by 
	\[
	f \cdot g(s_1, s_2, \ldots, s_n) = 
	\begin{cases} 
		f(2s_1, s_2, \ldots, s_n) & \text{if } s_1 \in [0, \frac{1}{2}] \\
		g(2s_1 - 1, s_2, \ldots, s_n) & \text{if } s_1 \in [\frac{1}{2}, 1]
	\end{cases}
	\]
	Specifically, $\pi_2(X, x_0)$ is called the second homotopy group.
\end{definition}
\begin{definition}[Homotopy Group]
	Let $X$ be a topological space, a subspace $A \subseteq X$ is a strong deformation retract if there exists a homotopy $h : X \times I \to X$ such that
\[
\begin{aligned}
	h(x, 0) &= x \quad \forall x \in X \\
	h(x, t) &= x \quad \forall x \in A, \forall t \in I \\
	h(x, 1) &\in A \quad \forall x \in X
\end{aligned}
\]
\end{definition}
	\section{Proof of Theorem}
	
	\textbf{Theorem 1:} There is no computer algorithm that can determine whether the second homotopy group of arbitrary constructive compact topological space is trivial or not.
	
	The sketch of the proof is shown in the following flowchart.

\begin{center}
	\includegraphics[scale=0.9]{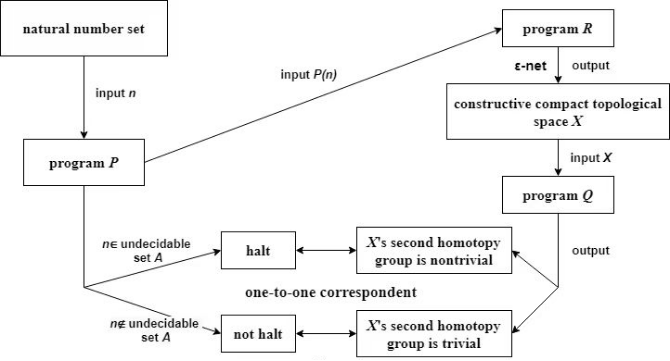}
\end{center}
	The proof of the theorem requires the following lemmas:
	
	\textbf{Lemma 1:} Let $X$ be a topological space, $A$ is the strong deformation retract of $X$. We have $\pi_n(X, a_0) = \pi_n(A, a_0) \quad \forall a_0 \in A, \forall n \geq 0$. This is a classical result in Algebraic Topology \cite{hatcher}.
	
	\textbf{Lemma 2:} $\pi_2(S^2) \cong \pi_2(\partial I^3) \cong \mathbb{Z}$, so they are nontrivial. Moreover, the second homotopy group for a cube is trivial.
	
	\textbf{Lemma 3:} Let $X$ be a constructive compact topological space, $\{\varepsilon_n\}_{n=0}^\infty$ is a sequence such that $\varepsilon_n = \frac{1}{2^n}$. For every $\varepsilon_n$, there exists a finite point set $S_{\varepsilon_n}$ which is a $\varepsilon$-net of $X$ due to the compactness of $X$. Then the union of these finite sets of points $S = \bigcup_{n=0}^\infty S_{\varepsilon_n}$ is dense in $X$. In other words, $S = X$.
	
	\textbf{Lemma 4:} A cube with a single central hollow has nontrivial second homotopy group. The proof of this Lemma is straightforward.
	
	\textbf{Lemma 5:} There does exist a program whose domain is undecidable \cite{shen}.
	
\subsection*{Proof of Theorem 1}
	We assume there exists an algorithm that can determine whether the second homotopy group of any given constructive compact topological space is trivial. We denote this algorithm as $Q$. $Q$ takes a constructive compact topological space as input and outputs "No" if its second homotopy group is trivial, otherwise it outputs "Yes".
	
	To derive a contradiction, we need to construct a program $P$ where the domain of $P$ is an enumerable and undecidable set. From Lemma 5, such program $P$ does exist. For each given positive integer $n$, the algorithm $P$ either does not halt at $n$ or halts at $n$ within a specific integer time. 
	
	We also need to construct an algorithm $R$ which takes ‘whether $P(n)$ halts’ as input and constructs a topological space by building $\varepsilon$-net on the cube $A = [0,1]^3$ using the method that depends on whether $P(n)$ halts or not. Specifically, $R$ has two methods for constructing the $\varepsilon$-net which we refer to as Method One and Method Two. If $P(n)$ does not halt, then $R$ uses Method One to construct the $\varepsilon$-net; otherwise, it uses Method Two (which goes as Method 1 till Halting moment of $P$ and then switches). After the construction of $\varepsilon$-net, we shall get a new compact constructive topological space. Note that $A$ is a closed and bounded set in Euclidean space, which implies $A$ is a compact metric space. Now we define a sequence $\{\varepsilon_m\}_{m=0}^\infty$ such that $\varepsilon_m = \frac{1}{2^m}$. 
	
	We construct a decision algorithm $F$ of $P$’s domain using the algorithms $P$, $R$, and $Q$ to obtain the contradiction by the following process:
	Firstly, given an arbitrary positive integer $n$ as the input for $P$ (which is also the input of $F$). $P$ runs for the first second. If $P$ does not halt in the first second, $P$ pauses and the result that $P$ does not halt with the input $n$ will be used as the input for $R$. Then $R$ begins to run for one second. During this second, $R$ handles the case $m=0$, generating an $\varepsilon$-net with $\varepsilon = 1$ using Method One. After that, $P$ runs for the second second with the initial input $n$. If $P$ has not halted in the second second, $P$ pauses again and the result that $P$ does not halt with the input $n$ will be used as the input for $R$ again. Then $R$ begins to run for the second second. During this second, $R$ handles the case $m=1$, generating an $\varepsilon$-net with $\varepsilon = \frac{1}{2}$ using Method One. 
	Now this process can be iterated if $P$ does not halt on the given $n$, increasing the value of $m$ with a step size of one. At this point, $R$ will eventually output $A$. Then we use $A$ as the input for $Q$ and run algorithm $Q$. In this case, $Q$ will output "No", which is the final output of $F$.
	
	Suppose that $P$ halts at a certain second with the input $n$. Clearly, the behavior of $R$ before $P$ halts has been described above. After $P$ halts, the result that $P$ does halt with the input $n$ will be used as the input for $R$. In this case, $R$ will handle the corresponding value of $m$, generating the corresponding $\varepsilon_m$-net using Method Two. Subsequently, $R$ will continue to run indefinitely with the value of $m$ gradually increasing (with a step size of 1), generating the corresponding $\varepsilon_m$-net using Method Two.  
	At this point, $R$ will eventually output a constructive compact topological space whose second homotopy group is nontrivial. Then using this output as the input for $Q$ and running algorithm $Q$, $Q$ will output "Yes", which is the final output of $F$ (The specific details of this explanation will be provided later in the text).
	
	Now we will introduce how Method One and Method Two specifically construct the $\varepsilon$-net and the constructive compact topological space. And we will divide it into the following two cases:

	(Method One) If $P$ does not halt on the positive integer $n$: 
	Let $B_m$ denote the set of points forming the $\varepsilon_m$-net. If $m=0$, we take $B_0 = \{0,1\}^3$. Having constructed $B_0, B_1, B_2, \ldots, B_{m-1}$, $B_m$ is the set 
	$$B_m = \left\{(x_m, y_m, z_m) \mid (x_m, y_m, z_m) \in \left\{k \cdot 2^{-m} \mid k = 0, 1, 2, \ldots, 2^m \right\} \right\} \backslash \bigcup_{j=0}^{m-1} B_j$$
	where $i = 1, 2, \ldots, 2^{m+1}-1$, $m \geq 1$. It is not difficult to show that $\text{card}(B_m) = 2^{m+1}-1$, $m \geq 1$ and $B_m$ is an $\varepsilon_m$-net of $A$. Since $P$ does not halt on the positive integer $n$, the value of $m$ will start from zero and continue to infinity without altering the construction method of $B_n$. Thus we will obtain $B_i (i = 0, 1, 2, \ldots)$. Then the output of program $R$ is $\lim_{n \to \infty}\bigcup_{i=1}^{n}B_n$ which is the whole space $[0,1]^3$ according to the Lemma 3.
	The second homotopy group of $A$ is trivial by Lemma 1 and Lemma 2. After $R$ outputs this result, we use it as the input for $Q$ and start running $Q$. The result of $Q$ 's execution will be "No".
	
	(Method Two) If $P$ halts on $n$ within a specific time: Clearly, before $P$ halts, $R$ constructs each $\varepsilon_m$-net using Method One. Therefore, we only need to consider the construction method of $R$ after $P$ stops. Suppose $B_m$ always denotes the $\varepsilon_m$-net generated by Method One and $B_m'$ always denotes the $\varepsilon_m$-net generated by Method Two. Now suppose that $B_0, B_1, B_2, \ldots, B_{m-1}$ $(m \geq 1)$ have been constructed using Method One and $P$ halts with input $n$ within the second before program $R$ starts to construct $B_m'$. Suppose $C$ is the center point of the cube $A$. Let                                                                              
$$d = \frac{1}{2} \min\{d(x_1, C), d(x_2, C), \ldots, d(x_k, C) \}, x_j \in \bigcup_{i=0}^{m-1} B_i, j = 1, 2, \ldots, k, k = \sum_{i=0}^{m-1}\text{card}(B_i).$$
	Then 
	$$B_i' = B_i \setminus N_d(C) \quad (i \geq m).$$
	(If $P$ halts in the first second, we can directly choose $d = \frac{1}{4}$). This is the method by which $R$ constructs the $\varepsilon$-net after $P$ halts. Clearly, $B_0, B_1, \ldots, B_{m-1}, B_m', B_{m+1}', \ldots$ is an $\varepsilon$-net of the space $A \setminus N_d(C)$. Then the closure of the union of these previous sets is the space $A \setminus N_d(C)$ if $P$ halts in the first or second second; is the space $A \setminus N_d(C) \cup \{\frac{1}{2}, \frac{1}{2}, \frac{1}{2}\}$ if $P$ halts after the third second. These two compact spaces are all the possible output results of algorithm $R$ when $P$ halts (The specific space depends on the exact second when $P$ halts). By Lemma 2 and Lemma 4, we can know that the second homotopy groups of these two spaces are non-trivial. Similar to before, after $R$ outputs this result, we use it as the input for $Q$ and start running $Q$. The result of $Q$ 's execution will be "Yes". 
	
	Now we have clearly constructed algorithm $F$. For each positive integer $n$ as the input for $P$:
	\begin{itemize}
		\item If $P(n)$ halts, then $F(n)$ outputs "Yes" which means $n$ is in the domain of $P$.
		\item If $P(n)$ never halts, then $F(n)$ outputs "No" which means $n$ is not in the domain of $P$.
	\end{itemize}
	In this case, $F(n)$ becomes a decision program of $P$ which contradicts the fact that the domain of $P$ is undecidable.
\qed                                                                                                                                                                                                                                                                                                                                                                                                                                                                                                                                     
	
	\section{Proof of Lemmas}
	
	\textbf{Proof of Lemma 2:} $\pi_2(S^2) \cong \mathbb{Z}$ is a direct corollary of a famous result that $\pi_n(S^n) \cong \mathbb{Z}$ \cite{hatcher}.
	
	\textbf{Proof of Lemma 3:} Given any point $x \in X$, if $x \in S$, there is nothing to prove; if $x \in X \setminus S$, given any $\varepsilon > 0$, there exists an $\varepsilon_n$ such that $0 < \varepsilon_n < \varepsilon$ and a point $y \in S_{\varepsilon_n} \subseteq S$ such that $d(x, y) < \varepsilon_n < \varepsilon$ which implies $y \in N_\varepsilon(x)$. Then $x$ is a limit point of $S$. Finally, $X \subseteq S$ follows. Since $S \subseteq X$ always holds, we proved that $S = X$.\qed
	
	\textbf{Proof of Lemma 4:} This is true because by strong deformation retraction, we can prove that this space is homotopy equivalent to $I^3$ with the hole inside. Using Lemma 2, we can obtain its second homotopy groups is nontrivial.\qed
	
	\section{Further Conclusions}
	Since the key Lemma 2 used in the proof can be generalized to the $n$-dimensional case, and for higher-dimensional cases, we can still use a similar $\varepsilon$-net construction to obtain the program $R$ and decision algorithm $F$, we can assert the following corollary:
	
	\textbf{Corollary 1:} There does not exist a computer algorithm that can determine if the homotopy group $\pi_n$ of any compact topological space is trivial, regardless of the order of the homotopy group.

	\section{Acknowledgement}
	We would like to express our sincere gratitude to Professor Vladimir Chernov, whose significant guidance and assistance were crucial in overcoming key challenges in this project, and to him and Viktor Chernov who formulated this project. We also extend our appreciation to the Neoscholar Company for their support. Special thanks to our teaching assistant Evans Huang for providing detailed explanations and helping us better organize and understand our study. Finally, we acknowledge the collective effort of our entire team whose dedication and hard work made the completion of this project possible.
	
\end{document}